%% file: sos.tex
\documentclass[11pt]{amsart} 
\usepackage{amsmath, amssymb, mathrsfs,pstricks}%, booktabs} 
\usepackage{url} % formats url's nicely using \url{}
\usepackage{multirow}
\usepackage{amsaddr}
\usepackage{fancyhdr}
{
\newtheorem{theorem}{Theorem}[section]
\newtheorem{lemma}[theorem]{Lemma}

\newtheorem{definition}[theorem]{Definition}
\newtheorem{proposition}[theorem]{Proposition}
\newtheorem{conjecture}[theorem]{Conjecture}
}
{\theoremstyle{remark}

}

\newcommand\lref[1]{Lemma~\ref{lem:#1}}
\newcommand\tref[1]{Theorem~\ref{thm:#1}}
\newcommand\cref[1]{Corollary~\ref{cor:#1}}
\newcommand\sref[1]{Section~\ref{sec:#1}}

\newcommand\pref[1]{Proposition~\ref{prop:#1}}

\newcommand\conref[1]{Conjecture~\ref{con:#1}}
\newcommand\fref[1]{Figure~\ref{fig:#1}}

\renewcommand\proof{\noindent\textsl{Proof. }}
\newcommand\sqr[2]{{\vbox{\hrule height.#2pt
    \hbox{\vrule width.#2pt height#1pt \kern#1pt
        \vrule width.#2pt}\hrule height.#2pt}}}
\renewcommand\qed{%
	\ifmmode\eqno\sqr53
	\else\nolinebreak\ \hfill\sqr53\medbreak\fi}

\numberwithin{equation}{section}

\newcommand\diff{\mathbin{\mkern-1.5mu\setminus\mkern-1.5mu}}
\newcommand\ZBIS{\bZ~$\in \{$\BIS$\}$}
\newcommand\rten{R_{10}}
\newcommand\mkttd{(M(K_{3,3}))^*}
\newcommand\mkfd{(M(K_5))^*}
\newcommand\mkttdp{(M(K_{3,3}^+))^*}
\newcommand\forma{{\Large \wp}}
\newcommand\ag{\ms {AG}(3,2)}
\newcommand\pef{\pu ef}
\newcommand\pu{\cup}

\newcommand\gef{G_e^f}
\newcommand\gfe{G_f^e}
\newcommand\Gef{G_{ef}}
\newcommand\GEF{G^{ef}}

\newcommand\vy{\mathbf y}
\newcommand\bI{{\bfseries \textsf{I}}}
\newcommand\bB{{\bfseries \textsf{B}}}
\newcommand\bS{{\bfseries \textsf{S}}}
\newcommand\bZ{{\bfseries \textsf{Z}}}
\newcommand\BIS{{\bfseries \textsf{B, I, S}}}

\newcommand\mc{\mathcal}
\newcommand\ms{\mathscr}
\newcommand\rr[1]{(\ref{#1})}
\newcommand\lp{\left(}
\newcommand\rp{\right)}

\newcommand\ass{:=}
\newcommand\ssum[1]{\sum_{ \substack{#1}}}
\newcommand\rref[1]{(\ref{#1})}
\newcommand\cl{\text{cl}}
\newcommand\set[1]{\left\{#1\right\}}

\newcommand\seymour{\cite{seymour-regular}}
\newcommand\fkg{\cite{fkg}}
\newcommand\pemantle{\cite{pemantle}}

\newcommand\sw{\cite{semple-welsh}}
\newcommand\wc{\cite{wagner-choe}}

\newcommand\gw{\cite{grimmett-winkler}}
\newcommand\ncrv{\cite{ncrv}}

\newcommand\kir{\cite{kirchhoff}}

\newcommand\kahnforests{\cite{kahnforests}}
\newcommand\kn{\cite{kahn-neiman}}

\newcommand\merino{\cite{merino}}
\newcommand\sflows{\cite{sflows}}

\newcommand\cocks{\cite{cocks}}
\newcommand\masters{\cite{Erickson-2008}}

\title[Negative correlation for spanning forests of graphs]{Sums of
  squares and negative correlation for spanning forests of series
  parallel graphs}

\author{Alejandro Erickson} \address{University of Victoria\\ate@uvic.ca}

\begin{document}
\maketitle

\input{sostext}

%http://amath.colorado.edu/documentation/LaTeX/reference/faq/bibstyles.html
\bibliographystyle{plain}%{mwnbib}
\bibliography{sos}
\nocite{wc}

\end{document}

%% file: sostext.tex
\begin{center}\textbf{Abstract}\\
  \vspace*{.5cm}
\end{center}

We provide new evidence that spanning forests of graphs satisfy the
same negative correlation properties as spanning trees, derived from
Lord Rayleigh's monotonicity property for electrical networks.  The
main result of this paper is that the Rayleigh difference for the
spanning forest generating polynomial of a series parallel graph can
be expressed as a certain positive sum of monomials times squares of
polynomials.  We also show that every regular matroid is
independent-set-Rayleigh if and only if every basis-Rayleigh binary
matroid is also independent-set-Rayleigh.

\section{Introduction}

The well-known theorem of Fortuin, Ginibre, and Kasteleyn in
\fkg~gives a locally verifiable sufficient condition for identifying
positively associated measures.  Unlike its well explored counterpart,
no such condition is known for negatively associated measures, as
defined by Kahn in \kn, nor is there one for the following special
case.  Consider a measure on a collection of sets so that the
probability of an event $\mc A$ occurring is $P(\mc A)$.  We say the
measure is \emph{negatively correlated} if
\begin{align}
  \label{matroidmeasure} P(\set{X:e,f\in X}) \leq P(\set{X:e\in
    X})P(\set{X:f\in X})
\end{align}
for every pair of distinct elements $e$ and $f$.

We are concerned with a certain family of measures that are positive
for bases, independent sets and spanning sets of matroids, represented
here by the letters \bB, \bI~ and \bS, respectively.  Let $\mc M$ be a
matroid and let $(y_g:g\in E)$ be a weighting of its ground set so
that the polynomial
\begin{align*}
  Z = \sum_{X} \vy^X,
\end{align*}
is a sum over \bZ-sets where \ZBIS~and $\vy^X = \lp \prod_{x\in X}
y_x\rp$.

Let $Z_e$ indicate the partial derivative $\frac{\partial Z}{\partial
  y_e}$.  For a given positive evaluation of the $y_g$s, suppose the
term $\vy^X$ is selected with probability $P(\set{X}) :=
\frac{\vy^X}{Z}$.  The monomials $y_eZ_e$ are precisely those of $Z$
which contain $y_e$, so that $\frac{y_eZ_e}{Z} = P(\set{X:e\in X})$.
The difference
\begin{align}
  \Delta Z \set{e,f} := Z_eZ_f-ZZ_{ef}
\end{align}
is called the \emph{\bZ-Rayleigh difference} and it is non-negative
for every positive evaluation of the $y_g$s if and only if
(\ref{matroidmeasure}) holds for the corresponding measures.  If
$\Delta Z\set{e,f} \ge 0$ for every pair of distinct edges $e$ and $f$
and every positive evaluation of the $y_g$s, then $\mc M$ is
\bZ-Rayleigh.

% [DO I NEED THIS?]A matroid is \bZ-negatively correlated if $\Delta Z
% \set{e,f}$ is non-negative when the $y_g$s are all evaluated at 1.

Graphs are \bB-Rayleigh as a result of Kirchhoff's laws for electrical
resistor networks (\kir) and the intuitive property, due to Lord Rayleigh, that
increasing the conductance of any resistor in the circuit does not
decrease the conductance between any two nodes.  Most of our efforts
are spent on the spanning forest analogue first circulated by Kahn in
the early 1990s (\kn).
\begin{conjecture}(\pemantle, \gw, \sw, \ncrv,
  \kahnforests)\label{con:giray} Graphs are \bI-Rayleigh.
\end{conjecture}
Independent efforts by Cocks in \cocks, Semple and Welsh in \sw~and
work by Wagner, especially \ncrv, prove that two-sums of \bI-Rayleigh
graphs are \bI-Rayleigh.  Grimmett and Winkler show in \gw~that graphs
on at most eight vertices and nine vertices with at most 18 edges have
a non-negative \bI-Rayleigh difference when the variables are
evaluated at $1$.  Cocks (\cocks) and Erickson (\masters) prove
independently that if all graphs satisfy this last condition, then
they are all \bI-Rayleigh as well.

Let $\mc G :=(V,E)$ be a graph with distinct edges $e$ and $f$ and let
$G$ denote the generating polynomial for edge sets of its spanning
forests, the \bI-sets.  Wagner conjectures that the \bI-Rayleigh
difference for $\mc G$ has the form
\begin{align}
  \label{spform}\Delta G\set{e,f} = \sum_{S\subseteq E} \vy^S
  A(S)^2,
\end{align}
where the sum is over sets $S$ which are contained in cycles through
both $e$ and $f$. For each $S$, the polynomial $A(S)$, equal to $
\sum_A c(S,e,f,C)\vy^{A-S}$, is a sum over all spanning forests $A$
such that $A\cup \set{e,f}$ contains a unique cycle $C$ which contains
$S$.  The signs $c(S,e,f,C)$, however, are not known.  The main result
is that the \bI-Rayleigh difference for any series-parallel graph can
be expressed this way.
% \label{thm:spsos} SP graphs are SOS
\begin{theorem}\label{thm:spsos}
  If $\mc G$ is a series-parallel graph, then $\mc G$ satisfies
  \rref{spform} for some choice of signs $c(S,e,f,C)$.
\end{theorem}
% make sure this is somewhere A computer program and some
% educated guesses for the signs shows that Wagner's conjecture
% holds for graphs on at most seven vertices (\masters).
In \sref{minordirect} we prove that if \rref{spform} holds for graphs
$\mc G$ and $\mc H$, then it holds for minors and direct sums of these
and in \sref{twosums} we present evidence that it also holds by taking two-sums.

Regular matroids are closely related to graphic matroids through
decomposition. In \sref{binreg} we prove the following relationship
between regular and binary matroids.
\begin{theorem}\label{thm:binregiray}
  The following are equivalent.
  \begin{enumerate}
  \item[(i)] Regular matroids are \bI-Rayleigh.
  \item[(ii)] Every \bB-Rayleigh binary matroid is also \bI-Rayleigh.
  \end{enumerate}
\end{theorem}

\section{Conjecture: Graphs are \bI-Rayleigh}\label{sec:minordirect}

Let $\mc G:=(V,E)$ be a graph whose spanning forests are denoted $\ms
F(\mc G)$.  More precisely, $\ms F(\mc G)$ is the collection of
acyclic subsets of $E$.  Their generating polynomial is
\begin{align*}
  F(\mc G;\vy) := \sum_{X\in \ms F(\mc G)} \vy^X,
\end{align*}
where $\vy := (y_g: g\in E)$ are indeterminates.  Write $G := F(\mc
G;\vy)$ and $H:=F(\mc H;\vy)$ and $K:=F(\mc K;\vy)$.  Braces and
commas are dropped from small sets of elements, as in $efg$ instead of
$\set{e,f,g}$.  We define the notation used in \conref{sos}.

\begin{definition}[S-sets, A-sets]\label{def:SOSnotation}
  Let $\mc G$ be a graph.  Let $\ms S$ be the collection of those sets
  $S \subseteq E-ef$ such that $S\pef$ is contained in some cycle of
  $\mc G$.  For each $S$ in $\ms S$, let $\ms A(S)$ be the collection
  of those spanning forests $A$ such that $A\subseteq E-ef$ and $S\pef
  \subseteq C\subseteq A\pef$ for a unique cycle, $C$. 
\end{definition}

Use a subscripted $G$ wherever the graph $\mc G$ needs to be
specified, as in $\ms A_G(S)$.  We refer to the elements of $\ms S$
and $\ms A(S)$ as S-sets and A-sets, respectively.  Throughout the
rest of this paper, given an S-set $S$ and one of its A-sets $A$, the
cycle $C$ is the unique cycle described in the above definition unless
otherwise noted.  The signs $c(S,e,f,C)$, used below, are written
$c(S,C)$ when $e$ and $f$ are understood.

% CONJECTURE
\begin{conjecture}[Wagner (private communication), Sum of
  Squares] \label{con:sos} Let $\mc G$ be a graph with distinct edges
  $e$ and $f$.  Then for some choice of signs $c(S,C) = \pm 1$,
  \begin{align}\label{eq:soscon}
    \Delta G\{e,f\} =\sum_{S\in \ms S} \vy^S \left( \sum_{A\in \ms
        A(S)} c(S,C)\vy^{A- S}\right)^2.
  \end{align}
\end{conjecture}
When $\mc G$ and $e$ and $f$ satisfy the above we say $\Delta G\set{e,f}$ is
SOS.  If $\mc G$ satisfies the above for every pair of distinct edges
$e$ and $f$ we say $\mc G$ is SOS.

\conref{sos} holds for the complete graph $K_7$, the cube and the
M\"obius ladder on eight vertices (Wagner, personal communication).
Other similarly sized graphs for which correct signs have not yet been
found, exhibit discrepancies on the order of tens out of tens of
thousands of terms.

Recall that $G_g$ is the partial derivative $\frac{\partial
  G}{\partial y_g}$.  When $g$ is not a loop, $G_g$ describes the
spanning forests of $\mc G$ contract $g$.  The analogue for deletion,
denoted $G^g$, is the evaluation at $y_g=0$.  We disclaim certain
omissions by remarking that loops have no more than a trivial effect
on our discussion of spanning forests and for the same reason we are
not concerned with whether or not $\mc G$ is connected.

\begin{lemma}\label{lem:sosminors}
  Let $\mc G$ be a graph with distinct edges $e,f$ and $g$.  If
  $\Delta G\{e,f\}$ is SOS then so are $\Delta G^g\{e,f\}$ and $\Delta
  G_g\{e,f\}$.
\end{lemma}
\proof From Section 4.4 of \ncrv, $\Delta G^g\{e,f\} = \lim_{y_g
  \rightarrow 0} \Delta G\{e,f\} $.  To show that this satisfies the
sum-of-squares form for $\Delta G^g\set{e,f}$ use
% \label{FGsos}
\begin{align}
  \label{Fgsos} \lim_{y_g \longrightarrow 0} \Delta G\{e,f\} =
  \ssum{S\in \ms S_{G} \\ g\not\in S} \vy^S \left(\ssum{A \in \ms
      A_{G}(S) \\ g\not\in A} c_G(S,C) \vy^{A- S} \right)^2.
\end{align}
A cycle containing a set $X$ is called an {\em $X$-cycle}.  An S-set
of $\mc G\diff g$ is contained in an $ef$-cycle of $\mc G\diff g$.
Clearly $\ms S_{G\diff g} \subseteq \{S: S\in \ms S_{G}, g\not\in
S\}$, the set indexing the outer sum of \rr{Fgsos}.  On the other
hand, given a set $\tilde S$ in $\{S: S\in \ms S_{G}, g\not\in S\}-
\ms S_{G\diff g}$, there are no $ef$-cycles containing $\tilde S$ and
not $g$.  Thus, there are no A-sets for $\tilde S$ which do not
contain $g$ and the inner sum of \rr{Fgsos} for these is empty.
Therefore, together, the sets indexing the sums in \rr{Fgsos} are the
S-sets and A-sets of $\mc G\diff g$.

The proof for $\Delta G_g\set{e,f}$ is slightly trickier due to the
fact that when $g$ is contracted, two cycles may be created from one.
Using $\lim_{y_g \rightarrow \infty} y_g^{-2}\Delta G \{e,f\} = \Delta
G_g \{e,f\}$ from \ncrv, terms of $\Delta G\set{e,f}$ without $y_g^2$
disappear, so we are left with
% \label{Fgsos}
\begin{align}
  \label{FGsos} \lim_{y_g \longrightarrow \infty} y_g^{-2} \Delta
  G\{e,f\} = \ssum{S\in \ms S_{G} \\ g\not\in S} \vy^S \left(\ssum{A
      \in \ms A_{G}(S) \\ g\in A} c_G(S,C) \vy^{A- (S\pu g)} \right)^2.
\end{align}
Observe that $g$ is not a chord of $C$ because the cycle $C$ is unique
in $A\pu ef$.  Thus, if every cycle containing $S$ has $g$ as a chord,
there are no A-sets in $\ms A(S)$ containing $g$.  Therefore we are
summing over S-sets not containing $g$ for which there is a cycle $C$
containing $S$ and $g$ is not a chord of $C$.  This is equal to $\ms
S_{G/g}$.

It remains to be shown that for an S-set $S$ of $\ms S_{G/g}$, the
inner sum of \rr{FGsos} is indexed by the desired A-sets.  Let $A'$ be an
element of $\ms A_{G/g}(S)$ and let $A \ass A'\pu g$.  By definition
$A'\pu e$ and $A'\pu f$ are forests of $\mc G/g$ and therefore $A\pu e$
and $A\pu f$ are forests of $\mc G$.  Furthermore there is a unique cycle
$C$ such that $S\cup ef\subseteq C\subseteq A\pef$ containing $S$, so
$A \in \ms A_G(S)$ and $g\in A$.

Conversely, suppose $A\in \ms A_G(S)$ and $g\in A$.  Clearly $A- g\in
\ms A_{G/g}(S)$, since $g$ cannot be a chord of $C$.  Therefore the
S-sets and A-sets of $\mc G/g$ are exactly those sets indexed by
\rr{FGsos}.\qed
% $\ms S_{G/g}$ and $\ms A_{G/g}(S)$ for $S\in \ms S_{G/g}$

For graphs $\mc H$ and $\mc K$, let the direct sum be any graph whose
spanning forests are generated by $HK$.  The sum-of-squares form also
holds by taking direct sums.

\begin{proposition}\label{prop:SOSdirectsum}
  If $\mc H$ and $\mc K$ are SOS graphs and $\mc G$ is their direct
  sum, then $\mc G$ is SOS.
\end{proposition}
\proof If $e\in E(\mc H)$ and $f\in E(\mc K)$, then there are no
cycles through $e$ and $f$ and hence no S-sets.  In this case $\Delta
G\set{e,f} = 0$.  Since $\mc H$ and $\mc K$ are both SOS, the case
where $\set{e,f}\subseteq E(\mc K)$ is eliminated by symmetry.

Let $e$ and $f$ be distinct edges in $E(\mc H)$.  Since $G=HK$, it is
easy to show that $\Delta G\set{e,f} = K^2 \Delta H\set{e,f}$, so that
\begin{align}
  \Delta G\set{e,f} = \sum_{S\in \ms S_H} \vy^S \left( \sum_{A\in \ms
      A_H(S)} c_H(S,C)\vy^{A- S} K\right)^2.\label{DeltaefHdirect}
\end{align}
The S-sets of $\mc G$ are equal to those of $\mc H$, since an
$ef$-cycle of $\mc G$ cannot contain an edge of $\mc K$.  Let $\ms A$
and $\ms B$ be collections of subsets of a set $X$ and define $ \ms
A\vee \ms B := \set{A\cup B: A\in \ms A, B\in \ms B}.$ Now, if $S$ is
an S-set of $\mc G$, then $\ms A_G(S) = \ms A_H(S) \vee \ms F(\mc K)$,
as required.\qed

%%%%%%%%%%%%%%%%%%%%%%%%%%%%%%%%%%%%%%%%%%%%%
% SECTION SERIES PARALLEL GRAPHS
\section{Series-Parallel Graphs}\label{sec:spgraphs}

Let $\mc H$ and $\mc K$ be graphs.  The two-sum, defined in \seymour,
of $\mc H$ and $\mc K$ along a common edge $g$ is denoted $\mc
H\oplus_g \mc K$.  In general there are up to two, non-isomorphic ways
of two-summing along $g$.  In spite of this fact the spanning forests
of two-sums are unique, so we do not make this distinction.

Denote the complete graph on three vertices by $K_3$ and let a
superscript $^*$ indicate matroid dual.  The graph $(K_3)^*$ consists
of three mutually parallel edges.  Define a {\em parallel extension}
of $\mc G$ to be $\mc G \oplus_g (K_3)^*$ for some edge $g$.
Similarly $\mc G \oplus_g K_3$ is called a {\em series extension}.  A
graph $\mc H$ is a {\em series-parallel extension} of $\mc G$ if it
can be obtained by a sequence of series and parallel extensions,
starting with $\mc G$.
% \begin{definition}
A graph is called {\em series-parallel} if it is a minor of a
series-parallel extension of $K_3$ or $(K_3)^*$.
% \end{definition}

We set out to prove \tref{spsos}, that every series-parallel graph is
SOS.  By \lref{sosminors} we need not consider proper minors of
series-parallel extensions of $K_3$ or $(K_3)^*$.  Let $\mc G := \mc H
\oplus_g \mc K$ for SOS graphs $\mc H$ and $\mc K$ and let $e$ and $f$
be distinct edges of $\mc G$.  We prove that if $\mc H$ and $\mc K$
are SOS and if $(H^g-H_g)H_g$ and $(K^g-K_g)K_g$ each satisfy a
similar sum-of-squares identity, then $\mc G$ is SOS and
$(G^h-G_h)G_h$ satisfies the same identity for every edge $h\in E(\mc
G)$.  The above mentioned identity is proved for all series-parallel
graphs in \lref{phispsos}.

There are three cases with respect to the locations of $e$ and $f$.
Either $e\in E(\mc H)-g$ and $f\in E(\mc K)-g$ or they are both in
$\mc H$ or in $\mc K$.  The last case is omitted by symmetry.  The
first holds for two-sums without any assumptions on $(H^g-H_g)H_g$ or
$(K^g-K_g)K_g$.

% \label{lem:sosehfk} case where e\in H and f\in K
\begin{lemma}\label{lem:sosehfk}
  Let $\mc G= \mc H \oplus_g \mc K$.  If $e\in E(\mc H)-g$ and $f\in
  E(\mc K)-g$, then $\Delta G\set{e,f}$ is SOS.
\end{lemma}
\proof From Theorem 5.8 of \ncrv, $\Delta G\{e,f\}=\Delta
H\{e,g\}\Delta K\{g,f\}$.  Since $\mc H$ and $\mc K$ are SOS,
\begin{align}\label{hksos}
   \begin{split}
   & \Delta H\{e,g\}\Delta K\{g,f\} \\
    =& \ssum{ S_H\in \ms S_H \\ S_K \in \ms
    S_K } \vy^{S_H\cup S_K } \left( \ssum{A_H\in \ms A_H(S_H) \\
      A_K\in \ms A_K(S_K)} c_G(S_G,C_G) \vy^{(A_H\cup A_K)- (S_H\cup
      S_K)} \right)^2    
 \end{split}
\end{align}
where we set $c_G(S_G,C_G):= c_H(S_H,C_H)c_K(S_K,C_K)$.

Notice that a cycle is an $ef$-cycle of $\mc G$ if and only if it is
the symmetric difference of an $eg$-cycle in $\mc H$ and a $gf$-cycle
in $\mc K$.  It is straightforward to show that \rref{hksos} is the
sum of squares we are expecting by showing that the outer and inner
sums index the S-sets and A-sets of $\mc G$, respectively.  \qed

The proof of the case where $\set{e,f}\subseteq E(\mc H)-g$ reduces to
proving a sum-of-squares form for $(K^g-K_g)K_g$.  For any graph $\mc
G$ and an edge $e$ let
\begin{align*}
  \Phi G\{e\} := (G^e-G_e)G_e.
\end{align*}

The proof of the following lemma is straight forward and similar to
Section 4.4 and Theorem 5.8 of \ncrv.

%%%%%%%%%% VERIFIED
% \label{lem:phiminor}
\begin{lemma}\label{lem:phiminor}
  Let $\mc G$ be a graph with distinct edges $e$ and $f$.  Then
  \begin{align}
    \label{phiminor}\Phi G\{e\} = \Phi G^f\{e\} + y_f
    \Psi G\{e|f\} + y_f^2\Phi G_f\{e\},
  \end{align}
  where
  \begin{align}
    \label{psi}\Psi G\{e|f\} = G_f^eG_e^f + G^{ef}G_{ef} -
    2G_e^fG_{ef}.
  \end{align}
  If $\mc G = \mc H \oplus_g \mc K$ and $e\in E(\mc H)$, then by
  setting $y_g := K^g/K_g -1$,
  % \label{phi2sum}
  \begin{align}
    \label{phi2sum} \Phi G\{ e \} = (K_g)^2 \Phi H\{e\}.
  \end{align}

\end{lemma}

%%%%%%%%%%%% END VERIFIED

To express the sum-of-squares form for $\Phi G\{e\}$ we need some
notation similar to that defined for \conref{sos}.  Note, however,
that Q-sets are required to be non-empty, unlike S-sets.  The
significance of this becomes clear later.

\begin{definition}[Q-sets, B-sets] Let $\mc G$ be a
  graph.  Let $\ms Q$ be the collection of those sets $Q$ such that
  $\varnothing \subset Q \subseteq E-e$ and $Q\pu e$ is contained in a
  cycle of $\mc G$.  For each $Q$ in $\ms Q$ let $\ms B(Q)$ be the
  collection of those spanning forests $B$ such that $B\subseteq E-e$
  and $Q\pu e\subseteq D\subseteq B\pu e$ for a unique cycle, $D$.
\end{definition}

Use a subscripted $G$
wherever the graph $\mc G$ needs to be specified, as in $\ms B_G
(Q)$. We refer to elements of $\ms Q$ and $\ms B(Q)$ as Q-sets and
B-sets, respectively.  Given a Q-set $Q$ and one of its B-sets $B$,
the cycle $D$ is the unique cycle described above. To avoid ambiguity,
the qualification, SOS, becomes $\Delta$-SOS.  If a graph $\mc G$ and
an edge $e$ satisfy the conclusion of the following lemma we say $\Phi
G\set{e}$ is $\Phi$-SOS.  If $\Phi G\set{e}$ is $\Phi$-SOS for every
edge $e$, then $\mc G$ is $\Phi$-SOS.

% \label{lem:phispsos}
\begin{lemma}\label{lem:phispsos}
  Let $\mc G$ be a series-parallel graph with an edge $e$.  With the
  above notation
  \begin{align}\label{phispsos}
    \Phi G\{e\} = \ssum{Q\in \ms Q} \vy^Q \left( \ssum{B\in \ms B(Q)}
      d(Q,D) \vy^{B- Q}\right)^2
  \end{align}

  for some choice of signs $d(Q,D) = \pm 1$.
\end{lemma}
\proof Recall that series-parallel graphs are minors of series
parallel extensions of $K_3$ and $(K_3)^*$.

Let $E(K_3) := \set{e,f,g}$ so that $\Phi K_3\{e\} = (1 + y_f + y_g +
y_fy_g)(1 + y_f + y_g) - (1 + y_f + y_g)^2 =y_f(y_g)^2 + y_g(y_f)^2 +
y_fy_g$ and $ \Phi (K_3)^*\set{e} = (1 + y_f + y_g)(1) - (1)^2= y_f +
y_g$.  Thus both $K_3$ and $(K_3)^*$ are $\Phi$-SOS.  Furthermore,
small modifications of \lref{sosminors} and \pref{SOSdirectsum} serve
to prove that the $\Phi$-sum-of-squares form holds by taking minors
and direct sums.  Therefore, by induction, it is enough to show that
two-sums of $\Phi$-SOS graphs are $\Phi$-SOS.

Let $\mc H$ and $\mc K$ be $\Phi$-SOS graphs such that $e\in E(\mc
H)-g$.  By \lref{phiminor} we have $\Phi G \{e\} = (K_g)^2 \Phi
H\set{e}$ and by the inductive hypothesis,
\begin{align}
  \label{aaa2} (K_g)^2 \Phi
H\set{e} = (K_g)^2 \ssum{Q\in \ms Q_H} \vy^Q \left( \ssum{B\in \ms
      B_H(Q)} d_H(Q,D) \vy^{B- Q}\right)^2,
\end{align}
in which $y_g = K^g/K_g - 1$.

A Q-set of \rref{phispsos} is contained in $\mc H$ or it is not.
Table \ref{tab:phisets}, which is divided according to this, shows the
bijections between index sets of \rref{phispsos} and \rref{aaa2},
highlighting the way they factor over the two-sum.  The case where
$g\not\in Q_H$ uses the fact that $d_H(Q,D)$ does not depend on $B-D$,
so we are able to group some B-sets of $\mc H$ (see Table
\ref{tab:phisets}).  The case where $g\in Q_H$ gives
\begin{align*}
  %\label{notmul}
  \Phi K\{g\}
  \ssum{Q: g\in Q} \vy^{Q- g} \left( \ssum{B\in \ms
      B_H(Q)}d_H(Q,D)\vy^{B- Q} \right)^2
\end{align*}
and it corresponds to Q-sets of $\mc G$ with edges in both factors.
For this reason Q-sets cannot be empty.  See \fref{PhiKPhiH} and Table
\ref{tab:phisets}.\qed

\begin{table}[h]
  \centering
  \begin{tabular}{p{1.1cm}||l|l||p{4.3cm}}
    notes &\multicolumn{2}{|l||}{\rr{aaa2},~~~~$(K_g)^2 \Phi H\set{e}$}&\rr{phispsos},~~~~$\Phi G\set{e}$\\
    \cline{2-4}
    & $\mc H$-part & $\mc K$-part & $\mc G$-part\\
    \hline
    \hline
    Q-sets
    &$Q_H:$ $g\not\in Q_H$
    & none
    & $Q_G:$ $Q_G\cap E(\mc K) = \varnothing$\\
    \hline
    &
    &
    &below, $B_H \in \ms B_H(Q_G)$\\
    \hline
    B-sets    
    & $B_H-g:$ $g\in D  $
    & $K^g-K_g$
    &$\set{B_H- g:g\in D}$ $\vee \left(\ms F(\mc K\diff g)- \ms F(\mc K/g)\right)$\\
    \hline
    \multirow{2}{1cm}{group terms}
    & $B_H-g:$ $g\in B_H- D$
    &       $K^g-K_g$
    &\multirow{2}{3cm}{$\set{B_H:g\not\in B_H, g\not\in \cl(B_H)}$ $\vee \ms F(\mc K\diff g)$}\\
    \cline{2-3}
    &$B_H:g\not\in B_H,$ $g\not\in \cl(B_H)$
    &$K_g$
    &\\
    \hline

   &$B_H:g\not\in B,$ $g\in \cl(B_H)$
   & $K_g$
   &$\set{B_H:g\not\in B_H, g\in \cl(B_H)}$ $\vee \ms F(\mc K/g)$\\
   \hline
   \hline
   Q-sets
   &$Q_H -  g:g\in Q_H$
   &$\ms Q_K$
   &   $Q_G: Q_G\cap E(\mc K) \neq \varnothing$\\
   \hline
   B-sets
   &$B_H-g:$ $B_H\in \ms B_H(Q_H)$
   &$B_K(Q_K)$
   &   $B_G\in\ms B_G(Q_G) $\\
   \hline
   \hline
  \end{tabular}
  \caption{Explicit bijections on the index sets of \rref{phispsos}
    and \rref{aaa2}.  The $\mc K$-part column accounts for $y_g$ and
    the $(K_g)^2$ factor. Write $g\in \cl(X)$ if and only if $g$
    completes a cycle in the set of edges $X$  and recall the
    notation $\ms A\vee \ms B = \set{A\cup B: A\in \ms A, B\in \ms
      B}$.  Finally, sets have been labelled naturally so that $Q_H\in
    \ms Q_H$, et cetera.}
  \label{tab:phisets}
\end{table}

%%%%%%%%%%%%%%%%%%%%%%%
%%%%%%%%%%%%%%%%%%%%%%% END OF PHI SP PROOF
%%%%%%%%%%%%%%%%%%%%%%%
%%%%%%%%%%%%%%%%%%%%%%%

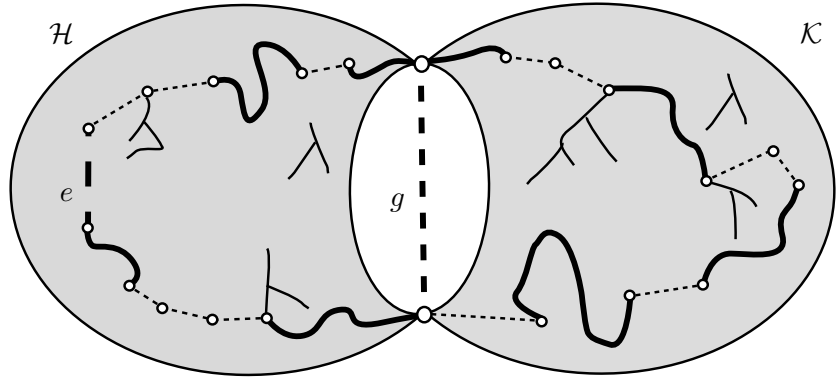
\begin{figure}[h] %  figure placement: here, top, bottom, or page
  \centering \input{phicase2}
  \caption[An element of $\ms B_G(Q_G)$.]{An element of $\ms B_G(Q_G)$
    may look like this.  Thin dashed lines represent edges of
    $Q_G=(Q_H- g) \cup Q_K$, thick solid lines complete a cycle $D_G$
    containing $Q_G \pu e$.}
  \label{fig:PhiKPhiH}
\end{figure}

We use \lref{phispsos} to prove in a similar way, that the
$\Delta$-sum-of-squares conjecture holds over two-sums when
$\set{e,f}\subseteq E(\mc H)$ and $\mc K$ is $\Phi$-SOS.

\begin{lemma}\label{lem:sosefh}
  Let $\mc G := \mc H \oplus_g \mc K$ and let $e$ and $f$ be distinct
  edges in $E(\mc H) -g$. If $\mc K$ is $\Phi$-SOS and $\mc H$ is
  $\Delta$-SOS, then $\Delta G\set{e,f}$ satisfies the
  $\Delta$-sum-of-squares form for some choice of signs $c_G(S,C)$.
\end{lemma}

\proof From Theorem 5.8 of \ncrv, $\Delta G\{e,f\} = (K_g)^2 \Delta H
\{e,f\},$ where $y_g = K^g/K_g - 1$.  By assumption we have
\begin{align}
  (K_g)^2 \Delta H \{e,f\} =(K_g)^2 \ssum{S \in \ms S_H} \vy^S \left(
    \ssum{A\in \ms A_H(S)} c_H(S,C)\vy^{A- S}
  \right)^2.\label{eq:KgDeltaHSOS}
\end{align}
%%%showing this is equal to \label{eq:soscon}

An S-set of \rref{eq:soscon} is contained in $\mc H$ or it is not.
Table \ref{tab:deltasets}, which is divided according to this, shows
the bijections between index sets of \rref{eq:soscon} and
\rr{eq:KgDeltaHSOS}, highlighting the way they factor over the
two-sum.  The case where $g\not\in S_H$ uses the fact that $c_H(S,C)$
does not depend on $A-C$, so we are able to group some A-sets of $\mc
H$ (see Table \ref{tab:deltasets}). 

The terms indexed by S-sets containing $g$ are
\begin{align}
  \label{notmuldelta}\Phi K\{g\} \ssum{S: g\in S} \vy^{S- g} \left(
    \ssum{A\in \ms A_H(S)}c_H(S,C)\vy^{A- S} \right)^2
\end{align}
and it corresponds to S-sets having edges in both $\mc H$ and $\mc K$.
The parts in $\mc K$ are the Q-sets of $\Phi K\set{g}$.  See Table
\ref{tab:deltasets} and again \fref{PhiKPhiH} noting this time that
$f$ is on the cycle in $\mc H$.\qed% \tag*{\sqr53}
\begin{table}[h]
  \centering
  \begin{tabular}{p{1.1cm}||l|l||p{4.3cm}}
    notes &\multicolumn{2}{|l||}{\rr{eq:KgDeltaHSOS},$~~~~ (K_g)^2 \Delta H\set{e,f}$}&\rr{eq:soscon},$~~~~ \Delta G\set{e,f}$\\
    \cline{2-4}
    & $\mc H$-part & $\mc K$-part & $\mc G$-part\\
    \hline
    \hline
    S-sets
    &$S_H:$ $g\not\in S_H$
    & none
    & $S_G: S_G\cap E(\mc K) = \varnothing$\\
    \hline
    &
    &
    &below, $A_H \in \ms A_H(S_G)$\\
    \hline
    A-sets
    & $A_H-g:$ $g\in C  $
    & $K^g-K_g$
    &$\set{A_H- g: g\in C}$ $\vee \left(\ms F(\mc
      K\diff g)- \ms F(\mc K/g) \right)$\\
\hline
    \multirow{2}{1cm}{group terms}
    & $A_H-g:$ $g\in A_H- C$
    &       $K^g-K_g$
    &\multirow{2}{3cm}{$\set{A_H:g\not\in A_H, g\not\in
        \cl(A_H)}$ $\vee \ms F(\mc K\diff g)$}\\
    \cline{2-3}
  &$A_H:g\not\in A_H,$ $g\not\in \cl(A_H)$
  &$K_g$
  &\\
\hline
   &$A_H:g\not\in A_H,$ $g\in \cl(A_H)$
   & $K_g$
   &$\set{A_H:g\not\in A_H, g\in \cl(A_H)}$ $\vee \ms F(\mc K/g)$\\
   \hline
   \hline
   S-sets
   &$ S_H -  g:$ $g\in S_H$
   &$\ms Q_K$
   &   $S_G: S_G\cap E(\mc K) \neq \varnothing$\\
   \hline
   A-sets
   &$A_H-g:$ $A_H\in \ms A_H(S_H)$
   &$\ms B_K(Q_K)$
   &   $A_G\in\ms A_G(S_G) $\\
   \hline
   \hline
  \end{tabular}
  \caption{
    Explicit bijections on the index sets of \rr{eq:soscon}
    and \rr{eq:KgDeltaHSOS}.   See notes at Table \ref{tab:phisets}.}
  \label{tab:deltasets}
\end{table}

We are finally in a position to prove \tref{spsos} which states that
series-parallel graphs are $\Delta$-SOS.

\proof(of \tref{spsos}) Let $\mc G$ be a series-parallel graph.
Either $\mc G$ is obtained by a sequence of series-parallel extensions
starting with $K_3$ or $(K_3)^*$, or $\mc G$ is a proper minor of one
of these.  By \lref{sosminors} we need only prove the theorem for the
first case.  It is straightforward to show that the base cases, $K_3$
and $(K_3)^*$ are $\Delta$-SOS.  Let $\mc G := \mc H \oplus_g \mc K$
where $\mc K$ is $K_3$ or $(K_3)^*$.  We assume that $\mc H$ is
$\Delta$-SOS.  If $e \in E(\mc H)$ and $f\in E(\mc K)$ then by
\lref{sosehfk}, $\Delta G\set{e,f}$ is $\Delta$-SOS.  If $\set{e,f}
\subseteq E(\mc H)$ or $\set{e,f}\subseteq E(\mc K)$ then
\lref{sosefh} is applicable, since $\mc H$ and $\mc K$ are $\Phi$-SOS
by \lref{phispsos}.  Thus, $\mc G$ is $\Delta$-SOS.\qed

% Implicit in the proofs of this section is some advice for choices of
% signs $c(S,C)$.  For example, consider the first inner term,
% $\ssum{A: g\in C} c_H(S,C) \vy^{A-(S\cup g)} (K^g-K_g) $, of
% \rr{splitSnoglast}.  In this case the S-set of $\mc H$ does not
% contain $g$ but the cycle in $\mc G$ is equal to $(C_H-g)\cup C_K$,
% for some $efg$-cycle of $\mc H$ and $g$-cycle of $\mc K$.  However,
% the signs for such terms are $c_G(S-g,(C_H-g)\cup C_K)$, which are
% equal to $c_H(S\cup g, C_H)$ for all choices of $C_K$.

%%%%%%%%%
% Two-Sums
%%%%%%%%%
\subsection{Two-sums of $\Delta$-SOS graphs}
\label{sec:twosums}
One might have hoped to prove, more generally, that if $\mc H$ and
$\mc K$ are $\Delta$-SOS, then $\mc H \oplus_g \mc K$ is as well.  The
problem lies in the fact that we are not assuming $\Phi H\set{g}$ and
$\Phi K\set{g}$ are $\Phi$-SOS.  To get around this we might try to
bootstrap this assumption by showing that it follows from the
induction hypothesis.  In fact, this looks promising and it is given
as the following conjecture.

\begin{conjecture}\label{con:myconjecture}%\label{con:myconjecture}
  If $\mc G$ is $\Delta$-SOS, then $\mc G$ is $\Phi$-SOS.
\end{conjecture}
Let $\mc G$ be a graph with distinct edges $e$ and $f$.  It is easy to
show that $\Delta G\set{e,f} = \gef\gfe - \Gef\GEF$ by using the fact
that $G = G^g + y_g G_g$ for any edge $g$.  Thus, recalling \lref{phiminor},
\begin{align}
  \notag \Psi G\set{e|f} = &G_f^eG_e^f + G^{ef}G_{ef} - 2G_e^fG_{ef} + (G^{ef}G_{ef}-G^{ef}G_{ef})\\
  =&\Delta G\set{e,f} + 2(G^{ef}G_{ef}-G_e^fG_{ef}).\label{newid}
\end{align}

We want to show that $\mc G$ is $\Phi$-SOS, whenever $\Phi G^f\set{e}$
and $\Phi G_f\set{e}$ are $\Phi$-SOS and $\Delta G\set{e,f}$ is
$\Delta$-SOS by showing that the expansion \rr{phiminor} can be
reduced to the desired $\Phi$-sum-of-squares form of $\Phi \mc
G\set{e}$, for some choice of signs $d_G(Q,D)$.

Dividing the $\Phi$-sum-of-squares sum for $\mc G$ into the two usual
cases where Q-sets do and do not contain $f$ yields
\begin{align}
  \label{Qnogsplit4}%\label{Qnogsplit4}
  \ssum{Q\in \ms Q_G\\ f\not\in Q} \vy^Q %LABEL%%%%%%
  \left( y_f\left(\ssum{B: f\in D} \forma(B) + \ssum{B:f\in B- D}
      \forma(B)\right)+ \ssum{B:f\not\in B\\f\not\in\cl(B)} \forma(B)
    + \ssum{B:f\not\in B\\f\in \cl(B)} \forma(B)
  \right)^2 \\
  \label{Qwithfdeg1}
  +y_f\left(\ssum{Q\in \ms Q\\f\in Q} \vy^{Q- f} \left( \ssum{B\in \ms
        B(Q)} d(Q,D) \vy^{B- Q}\right)^2\right),
\end{align}
where $\forma(B)$ stands for $ d(Q,D) \vy^{B- (Q\pu f)}$.  We are left
with comparing the coefficients of two polynomials in $y_f$, namely,
\rr{phiminor} and \rr{Qnogsplit4}-\rr{Qwithfdeg1}.  The degree $0$ and
$2$ terms come from \rr{Qnogsplit4} and the degree 1 terms are a
combination of \rr{Qwithfdeg1} with the cross terms of
\rr{Qnogsplit4}.

The methods in this section are easily adapted to showing that the
degree $0$ and $2$ terms are the sum-of-squares forms of $\Phi
G^f\set{e}$ and $\Phi G_f\set{e}$, respectively, and that
\rr{Qwithfdeg1} is the sum-of-squares form of $\Delta G\set{e,f}$,
which accounts for the first term in \rr{newid}.  We are left with
showing that the cross terms of \rr{Qnogsplit4} are equal to
$G^{ef}G_{ef}-G_e^fG_{ef}$.  The proof of the following proposition is
cumbersome and again similar to what we have seen.
	
\begin{proposition}
  If for each pair $e$ and $f$ of distinct edges,
  $G^{ef}G_{ef}-G_e^fG_{ef}$ is equal to
  \begin{align*}
    \ssum{Q\in \ms Q_G\\ f\not\in Q} \vy^Q %LABEL%%%%%%
    \left(\ssum{B: f\in D} \forma(B) + \ssum{B:f\in B- D}
      \forma(B)\right) \left(\ssum{B:f\not\in
        B\\f\not\in\cl(B)} \forma(B) + \ssum{B:f\not\in
        B\\f\in \cl(B)} \forma(B) \right),
  \end{align*}
  for a certain choice of signs, then \conref{myconjecture} is
  true.
\end{proposition}
	
\section{Binary matroids}\label{sec:binreg}
Graphic matroids are indeed interesting on their own, however, it is
worth being reminded of their role in the decomposition of regular
matroids.  A few relevant facts from \wc~and \seymour~are listed. For
undefined terms see \wc~and \seymour.
\begin{enumerate}
\item \label{binaryiffnos8} A binary matroid is \bB-Rayleigh if and
  only if it has no $S_8$ minor (\wc).
\item The affine geometry $\ag$ is a splitter for the class of binary matroids
  containing no $S_8$ minor (Seymour unpublished, Appendix D,
  \merino).
\item \label{sflows} If a binary matroid contains neither $\mc S_8$
  nor $\ag$ as a minor, then it can be constructed from direct sums
  and two-sums of regular matroids, the fano matroid, $F_7$, and its
  dual, $(F_7)^*$ (\sflows).
\item \label{subclass} From (3) it follows that a binary,
  three-connected matroid with no $S_8$ minor is regular or isomorphic
  to $F_7$, $(F_7)^*$ or $\ag$.
\item \label{r10r12} A three-connected regular matroid which is
  neither graphic nor co-graphic contains either $R_{10}$ or $R_{12}$
  as a minor (\seymour).
\item \label{seymour} Regular matroids decompose over direct sums,
  two-sums and three-sums into graphic and co-graphic matroids and
  $\rten$ (\seymour).
\end{enumerate}
	
We derive the following theorem.
\begin{theorem}\label{thm:binregiray}
  The following are equivalent.
  \begin{enumerate}
  \item[(i)] Regular matroids are \bI-Rayleigh.
  \item[(ii)] A binary matroid is \bI-Rayleigh if and only if it is
    \bB-Rayleigh.
  \end{enumerate}
\end{theorem}
\proof Assume (ii) and observe that by \rr{sflows}, regular matroids
are a subclass of binary matroids with no $S_8$ and no $\ag$ minor.
Since (ii) implies that binary matroids with no $S_8$ minor are
\bI-Rayleigh, (i) must be true.

Conversely, assume (i) and let $\mc M$ be a minor-minimal counter
example.  It is easy to show that if $\mc M$ is \bI-Rayleigh then it
must be \bB-Rayleigh, so we may assume that $\mc M$ is \bB-Rayleigh
and not \bI-Rayleigh.  By Theorem 5.8 of \ncrv~and its minimality,
$\mc M$ is three-connected.  Now from \rr{subclass}, $\mc M$ is either
$F_7$, $(F_7)^*$, $\ag$ or it is regular.  But $F_7$ and $(F_7)^*$ are
minors of $\ag$ which is \bI-Rayleigh (\sw~page 12).  Therefore $\mc
M$ must be regular, but that contradicts (i), so (i) implies (ii).\qed

\section{Concluding Remarks}
Whether or not regular matroids are \bI-Rayleigh is open, so we verify
the fact for small matroids.  In particular, regular matroids on up to
nine elements are \bI-Rayleigh.  Denote the graphic matroid of $\mc G$
by $M(\mc G)$.  A simple calculation in Maple shows that $\mkttd$ is
\bI-Rayleigh by first subtracting $\Delta I(\mkttd) \set{e,f} - \Delta
B(\mkttd)\set{e,f}$.  The resulting difference has four negative
terms, however, a small algebraic manipulation makes these disappear
into squares so that any positive evaluation of $\Delta I(\mkttd)
\set{e,f} - \Delta B(\mkttd)\set{e,f}$ is non-negative.  Thus $\Delta
I(\mkttd) \set{e,f} \ge 0$ for every positive weighting.  Let
$K_{3,3}^+$ be $K_{3,3}$ plus an edge not parallel to any others.
Three-connected regular matroids on at most 10 elements are either
$\mkttd$, $\mkfd$, $\mkttdp$, $\rten$ or graphic on at most six
vertices, which is the upper bound for such a graph.  Wagner has
verified that $K_6$ is $\Delta$-SOS which shows that regular matroids
on at most $9$ elements are \bI-Rayleigh.  Furthermore Cocks proves
that $\mkfd$ is \bI-Rayleigh (\cocks), so the only two obstructions to
showing (i) for 10 elements are $\mkttdp$ and $R_{10}$.  Unfortunately
for these last two, the method of subtracting $\Delta B\set{e,f}$ from
$\Delta I\set{e,f}$ yields not four, but tens of negative terms.

Semple and Welsh also ask whether graphs are \bS-Rayleigh.  This is
equivalent to co-graphic matroids being \bI-Rayleigh and it is
necessary for showing that regular matroids also possess the property.
Is there a sum-of-squares form for the \bS-Rayleigh difference,
analogous to that of \conref{sos}?

\paragraph{Acknowledgements} I thank David Wagner for his insight into
this problem and giving very helpful feedback on an early draft.

%%% Local Variables: 
%%% mode: latex
%%% TeX-master: "sos"
%%% End: 

%% file: phicase2.tex
%LaTeX with PSTricks extensions
%%Creator: inkscape 0.46
%%Please note this file requires PSTricks extensions
\psset{xunit=.5pt,yunit=.5pt,runit=.5pt}
\begin{pspicture}(634.95001221,297.67999268)
{
\newrgbcolor{curcolor}{0.86274511 0.86274511 0.86274511}
\pscustom[linestyle=none,fillstyle=solid,fillcolor=curcolor]
{
\newpath
\moveto(630.01837838,147.14182806)
\curveto(628.63055426,69.87114911)(547.75928605,8.59097675)(449.5021547,10.35572906)
\curveto(351.24502335,12.12048137)(272.62645516,76.2651792)(274.01427928,153.53585816)
\curveto(275.4021034,230.80653711)(356.27337161,292.08670947)(454.53050296,290.32195716)
\curveto(552.78763431,288.55720485)(631.4062025,224.41250702)(630.01837838,147.14182806)
\closepath
}
}
{
\newrgbcolor{curcolor}{0 0 0}
\pscustom[linewidth=1.87499993,linecolor=curcolor]
{
\newpath
\moveto(630.01837838,147.14182806)
\curveto(628.63055426,69.87114911)(547.75928605,8.59097675)(449.5021547,10.35572906)
\curveto(351.24502335,12.12048137)(272.62645516,76.2651792)(274.01427928,153.53585816)
\curveto(275.4021034,230.80653711)(356.27337161,292.08670947)(454.53050296,290.32195716)
\curveto(552.78763431,288.55720485)(631.4062025,224.41250702)(630.01837838,147.14182806)
\closepath
}
}
{
\newrgbcolor{curcolor}{0.8509804 0.8509804 0.8509804}
\pscustom[linestyle=none,fillstyle=solid,fillcolor=curcolor]
{
\newpath
\moveto(362.93977438,146.13960806)
\curveto(361.55195026,68.86892911)(280.68068205,7.58875675)(182.4235507,9.35350906)
\curveto(84.16641935,11.11826137)(5.54785116,75.2629592)(6.93567528,152.53363816)
\curveto(8.3234994,229.80431711)(89.19476761,291.08448947)(187.45189896,289.31973716)
\curveto(285.70903031,287.55498485)(364.3275985,223.41028702)(362.93977438,146.13960806)
\closepath
}
}
{
\newrgbcolor{curcolor}{0 0 0}
\pscustom[linewidth=1.87499993,linecolor=curcolor]
{
\newpath
\moveto(362.93977438,146.13960806)
\curveto(361.55195026,68.86892911)(280.68068205,7.58875675)(182.4235507,9.35350906)
\curveto(84.16641935,11.11826137)(5.54785116,75.2629592)(6.93567528,152.53363816)
\curveto(8.3234994,229.80431711)(89.19476761,291.08448947)(187.45189896,289.31973716)
\curveto(285.70903031,287.55498485)(364.3275985,223.41028702)(362.93977438,146.13960806)
\closepath
}
}
{
\newrgbcolor{curcolor}{1 1 1}
\pscustom[linestyle=none,fillstyle=solid,fillcolor=curcolor]
{
\newpath
\moveto(368.78377102,149.34848505)
\curveto(367.84310495,96.97448057)(343.51449846,54.89128558)(314.47884582,55.41278195)
\curveto(285.44319318,55.93427832)(262.64146495,98.86396017)(263.58213102,151.23796465)
\curveto(264.52279709,203.61196913)(288.85140358,245.69516412)(317.88705622,245.17366775)
\curveto(346.92270886,244.65217138)(369.72443709,201.72248953)(368.78377102,149.34848505)
\closepath
}
}
{
\newrgbcolor{curcolor}{0 0 0}
\pscustom[linewidth=1.87499994,linecolor=curcolor]
{
\newpath
\moveto(368.78377102,149.34848505)
\curveto(367.84310495,96.97448057)(343.51449846,54.89128558)(314.47884582,55.41278195)
\curveto(285.44319318,55.93427832)(262.64146495,98.86396017)(263.58213102,151.23796465)
\curveto(264.52279709,203.61196913)(288.85140358,245.69516412)(317.88705622,245.17366775)
\curveto(346.92270886,244.65217138)(369.72443709,201.72248953)(368.78377102,149.34848505)
\closepath
}
}
{
\newrgbcolor{curcolor}{0 0 0}
\pscustom[linewidth=4.5,linecolor=curcolor,linestyle=dashed,dash=13.5 13.5]
{
\newpath
\moveto(317.12638296,245.54773709)
\lineto(319.84258296,54.95597709)
\lineto(319.81283296,53.29962709)
}
}
{
\newrgbcolor{curcolor}{0 0 0}
\pscustom[linewidth=4.375,linecolor=curcolor]
{
\newpath
\moveto(263.66084296,246.87665709)
\curveto(263.66084296,246.87665709)(260.88776296,229.13216709)(277.60374296,237.42228709)
\curveto(294.31972296,245.71240709)(288.15269296,243.98237709)(295.55754296,246.30377709)
\curveto(302.96238296,248.62515709)(318.22021296,244.05593709)(318.22021296,244.05593709)
}
}
{
\newrgbcolor{curcolor}{0 0 0}
\pscustom[linewidth=4.375,linecolor=curcolor]
{
\newpath
\moveto(381.52139296,249.66857709)
\curveto(381.52139296,249.66857709)(378.66372296,261.37821709)(360.12958296,254.34795709)
\curveto(341.59544296,247.31769709)(318.87768296,246.49851709)(318.87768296,246.49851709)
}
}
{
\newrgbcolor{curcolor}{1 1 1}
\pscustom[linestyle=none,fillstyle=solid,fillcolor=curcolor]
{
\newpath
\moveto(323.24882822,244.60927722)
\curveto(323.19135284,241.40918719)(320.73886746,238.85522286)(317.7745305,238.90846399)
\curveto(314.81019354,238.96170512)(312.45100154,241.60208983)(312.50847692,244.80217987)
\curveto(312.5659523,248.0022699)(315.01843768,250.55623422)(317.98277464,250.50299309)
\curveto(320.9471116,250.44975196)(323.3063036,247.80936725)(323.24882822,244.60927722)
\closepath
}
}
{
\newrgbcolor{curcolor}{0 0 0}
\pscustom[linewidth=1.87499991,linecolor=curcolor]
{
\newpath
\moveto(323.24882822,244.60927722)
\curveto(323.19135284,241.40918719)(320.73886746,238.85522286)(317.7745305,238.90846399)
\curveto(314.81019354,238.96170512)(312.45100154,241.60208983)(312.50847692,244.80217987)
\curveto(312.5659523,248.0022699)(315.01843768,250.55623422)(317.98277464,250.50299309)
\curveto(320.9471116,250.44975196)(323.3063036,247.80936725)(323.24882822,244.60927722)
\closepath
}
}
{
\newrgbcolor{curcolor}{0 0 0}
\pscustom[linewidth=1.875,linecolor=curcolor,linestyle=dashed,dash=3.75 3.75]
{
\newpath
\moveto(408.87513296,50.24671709)
\lineto(377.08928296,52.33300709)
\lineto(319.56362296,55.63929709)
}
}
{
\newrgbcolor{curcolor}{0 0 0}
\pscustom[linewidth=4.375,linecolor=curcolor]
{
\newpath
\moveto(200.64682196,52.27132709)
\curveto(200.64682196,52.27132709)(211.40154296,36.12467709)(227.40498296,38.90522709)
\curveto(243.40841296,41.68577709)(238.76563296,56.49546709)(255.98485296,58.64058709)
\curveto(273.20405296,60.78569709)(255.07399296,42.08985709)(287.71629296,48.86672709)
\curveto(320.35858296,55.64360709)(320.33654296,54.41680709)(320.33654296,54.41680709)
}
}
{
\newrgbcolor{curcolor}{1 1 1}
\pscustom[linestyle=none,fillstyle=solid,fillcolor=curcolor]
{
\newpath
\moveto(325.21272942,54.85947232)
\curveto(325.15525404,51.65938229)(322.70276866,49.10541796)(319.7384317,49.15865909)
\curveto(316.77409474,49.21190022)(314.41490274,51.85228493)(314.47237812,55.05237497)
\curveto(314.5298535,58.252465)(316.98233888,60.80642932)(319.94667584,60.75318819)
\curveto(322.9110128,60.69994706)(325.2702048,58.05956235)(325.21272942,54.85947232)
\closepath
}
}
{
\newrgbcolor{curcolor}{0 0 0}
\pscustom[linewidth=1.87499991,linecolor=curcolor]
{
\newpath
\moveto(325.21272942,54.85947232)
\curveto(325.15525404,51.65938229)(322.70276866,49.10541796)(319.7384317,49.15865909)
\curveto(316.77409474,49.21190022)(314.41490274,51.85228493)(314.47237812,55.05237497)
\curveto(314.5298535,58.252465)(316.98233888,60.80642932)(319.94667584,60.75318819)
\curveto(322.9110128,60.69994706)(325.2702048,58.05956235)(325.21272942,54.85947232)
\closepath
}
}
{
\newrgbcolor{curcolor}{0 0 0}
\pscustom[linewidth=1.87499991,linecolor=curcolor]
{
\newpath
\moveto(461.73687423,223.89706268)
\curveto(459.14276382,224.41538039)(439.78007546,203.08871098)(428.98033706,194.66125045)
\curveto(417.50238566,185.70456373)(425.75404688,178.38128268)(433.10735812,168.07681332)
}
}
{
\newrgbcolor{curcolor}{0 0 0}
\pscustom[linewidth=1.87499991,linecolor=curcolor]
{
\newpath
\moveto(424.19884485,184.80576612)
\curveto(426.19194002,182.01274455)(408.55093595,163.61277953)(402.24214372,158.68979706)
\curveto(399.3708453,156.44921766)(399.79865475,156.1451962)(397.54988891,153.80339126)
}
}
{
\newrgbcolor{curcolor}{0 0 0}
\pscustom[linewidth=1.87499991,linecolor=curcolor]
{
\newpath
\moveto(442.96786427,204.35139931)
\curveto(441.03706675,202.84472903)(458.21394419,178.83984288)(460.87414839,175.8625819)
\curveto(462.37874167,174.17866147)(463.88332495,172.49474123)(465.38791841,170.8108308)
}
}
{
\newrgbcolor{curcolor}{0 0 0}
\pscustom[linewidth=1.87499991,linecolor=curcolor]
{
\newpath
\moveto(109.47719306,219.35948332)
\curveto(117.31803162,210.5840956)(101.2990002,191.3383174)(96.60751981,186.45270782)
}
}
{
\newrgbcolor{curcolor}{0 0 0}
\pscustom[linewidth=1.87499991,linecolor=curcolor]
{
\newpath
\moveto(107.61561984,201.16707014)
\curveto(109.76542967,198.76103673)(118.60483519,188.3443425)(121.09740116,182.69907182)
\curveto(125.82880182,177.40372331)(105.98030436,183.21028188)(101.06179163,178.08824115)
\curveto(99.68800975,176.65763044)(96.31962489,174.38774283)(94.83519435,173.22939299)
}
}
{
\newrgbcolor{curcolor}{0 0 0}
\pscustom[linewidth=1.87499991,linecolor=curcolor]
{
\newpath
\moveto(553.40947727,222.58671846)
\curveto(552.95893936,222.11754343)(541.89014015,205.54885221)(537.64962423,199.67323477)
\curveto(536.31403713,197.82267047)(534.52145109,196.4156343)(532.9573796,194.78683879)
}
}
{
\newrgbcolor{curcolor}{0 0 0}
\pscustom[linewidth=1.87499991,linecolor=curcolor]
{
\newpath
\moveto(550.00198091,237.39727662)
\curveto(548.81522104,235.75293298)(558.49730479,211.74161473)(561.71142658,205.70595366)
\curveto(563.28570273,202.74968678)(562.44149712,202.27415008)(563.12677264,199.05296487)
}
}
{
\newrgbcolor{curcolor}{0 0 0}
\pscustom[linewidth=1.87499991,linecolor=curcolor]
{
\newpath
\moveto(200.97488377,104.99819016)
\curveto(204.25153531,100.09628897)(200.30481484,67.69037888)(200.08240892,55.30735659)
}
}
{
\newrgbcolor{curcolor}{0 0 0}
\pscustom[linewidth=1.87499991,linecolor=curcolor]
{
\newpath
\moveto(200.46914805,76.8400528)
\curveto(200.4960205,76.79985376)(222.95510204,68.57055467)(227.84921918,63.09312481)
\curveto(229.12786678,61.6620626)(230.87820836,60.82951635)(232.39271839,59.69773685)
}
}
{
\newrgbcolor{curcolor}{0 0 0}
\pscustom[linewidth=1.87499991,linecolor=curcolor]
{
\newpath
\moveto(531.56936631,157.25619445)
\curveto(535.41252586,152.95496693)(564.65242941,151.27768637)(571.10512063,136.66336527)
}
}
{
\newrgbcolor{curcolor}{0 0 0}
\pscustom[linewidth=1.87499991,linecolor=curcolor]
{
\newpath
\moveto(551.36698266,148.61614283)
\curveto(554.83748227,144.73201465)(555.53188552,124.14031534)(555.3155184,112.0935188)
}
}
{
\newrgbcolor{curcolor}{0 0 0}
\pscustom[linewidth=4.5,linecolor=curcolor,linestyle=dashed,dash=13.5 13.5]
{
\newpath
\moveto(65.41822396,120.69847709)
\lineto(65.78316796,196.70418709)
}
}
{
\newrgbcolor{curcolor}{0 0 0}
\pscustom[linewidth=1.875,linecolor=curcolor,linestyle=dashed,dash=3.75 3.75]
{
\newpath
\moveto(65.23562896,196.59278709)
\lineto(110.42874986,223.81591709)
\lineto(160.57032996,231.25002709)
}
}
{
\newrgbcolor{curcolor}{1 1 1}
\pscustom[linestyle=none,fillstyle=solid,fillcolor=curcolor]
{
\newpath
\moveto(69.57875783,195.72593691)
\curveto(69.54286164,193.72732406)(67.85903786,192.13498062)(65.8202226,192.17159886)
\curveto(63.78140735,192.20821711)(62.15584984,193.85999887)(62.19174603,195.85861172)
\curveto(62.22764222,197.85722457)(63.911466,199.44956801)(65.95028126,199.41294976)
\curveto(67.98909652,199.37633152)(69.61465402,197.72454976)(69.57875783,195.72593691)
\closepath
}
}
{
\newrgbcolor{curcolor}{0 0 0}
\pscustom[linewidth=1.87500004,linecolor=curcolor]
{
\newpath
\moveto(69.57875783,195.72593691)
\curveto(69.54286164,193.72732406)(67.85903786,192.13498062)(65.8202226,192.17159886)
\curveto(63.78140735,192.20821711)(62.15584984,193.85999887)(62.19174603,195.85861172)
\curveto(62.22764222,197.85722457)(63.911466,199.44956801)(65.95028126,199.41294976)
\curveto(67.98909652,199.37633152)(69.61465402,197.72454976)(69.57875783,195.72593691)
\closepath
}
}
{
\newrgbcolor{curcolor}{0 0 0}
\pscustom[linewidth=4.375,linecolor=curcolor]
{
\newpath
\moveto(96.17440996,77.46435709)
\curveto(96.17440996,77.46435709)(109.14389856,82.14018709)(98.95851366,95.82223709)
\curveto(88.77311896,109.50427709)(74.58784496,105.46387709)(69.70269096,106.77880709)
\curveto(64.81754596,108.09374709)(64.43550096,120.98611709)(64.43550096,120.98611709)
}
}
{
\newrgbcolor{curcolor}{1 1 1}
\pscustom[linestyle=none,fillstyle=solid,fillcolor=curcolor]
{
\newpath
\moveto(68.90694783,120.63581691)
\curveto(68.87105164,118.63720406)(67.18722786,117.04486062)(65.1484126,117.08147886)
\curveto(63.10959735,117.11809711)(61.48403984,118.76987887)(61.51993603,120.76849172)
\curveto(61.55583222,122.76710457)(63.239656,124.35944801)(65.27847126,124.32282976)
\curveto(67.31728652,124.28621152)(68.94284402,122.63442976)(68.90694783,120.63581691)
\closepath
}
}
{
\newrgbcolor{curcolor}{0 0 0}
\pscustom[linewidth=1.87500004,linecolor=curcolor]
{
\newpath
\moveto(68.90694783,120.63581691)
\curveto(68.87105164,118.63720406)(67.18722786,117.04486062)(65.1484126,117.08147886)
\curveto(63.10959735,117.11809711)(61.48403984,118.76987887)(61.51993603,120.76849172)
\curveto(61.55583222,122.76710457)(63.239656,124.35944801)(65.27847126,124.32282976)
\curveto(67.31728652,124.28621152)(68.94284402,122.63442976)(68.90694783,120.63581691)
\closepath
}
}
{
\newrgbcolor{curcolor}{0 0 0}
\pscustom[linewidth=1.875,linecolor=curcolor,linestyle=dashed,dash=3.75 3.75]
{
\newpath
\moveto(95.67025296,77.08759709)
\lineto(122.63929096,59.93385709)
\lineto(161.11977796,50.90804709)
\lineto(202.06308496,52.44577709)
}
}
{
\newrgbcolor{curcolor}{0 0 0}
\pscustom[linewidth=1.875,linecolor=curcolor,linestyle=dashed,dash=3.75 3.75]
{
\newpath
\moveto(227.36228296,237.62738709)
\lineto(263.85610296,244.54891709)
}
}
{
\newrgbcolor{curcolor}{0 0 0}
\pscustom[linewidth=1.875,linecolor=curcolor,linestyle=dashed,dash=3.75 3.75]
{
\newpath
\moveto(382.91236296,249.98757709)
\lineto(418.44465296,245.56090709)
\lineto(459.75089296,225.11887709)
}
}
{
\newrgbcolor{curcolor}{0 0 0}
\pscustom[linewidth=1.875,linecolor=curcolor,linestyle=dashed,dash=3.75 3.75]
{
\newpath
\moveto(532.77056296,156.37228709)
\lineto(583.19784296,179.71291709)
\lineto(603.17294296,152.83472709)
}
}
{
\newrgbcolor{curcolor}{0 0 0}
\pscustom[linewidth=4.375,linecolor=curcolor]
{
\newpath
\moveto(408.65484296,52.21697709)
\curveto(408.65484296,52.21697709)(382.45500296,62.50506709)(389.99205296,72.18723709)
\curveto(397.52910296,81.86939709)(400.40133296,105.13444709)(404.13681296,108.13534709)
\curveto(407.87228296,111.13622709)(409.91078296,122.14434709)(421.43313296,114.57423709)
\curveto(432.95550296,107.00414709)(438.49456296,73.77047709)(439.03084296,69.46568709)
\curveto(439.56712296,65.16087709)(441.63512296,43.64787709)(444.03363296,40.53681709)
\curveto(446.43213296,37.42576709)(447.63690296,36.17692709)(460.45215296,32.26518709)
\curveto(473.26741296,28.35344709)(475.83472296,68.80466709)(475.83472296,68.80466709)
}
}
{
\newrgbcolor{curcolor}{0 0 0}
\pscustom[linewidth=4.375,linecolor=curcolor]
{
\newpath
\moveto(531.23884296,78.85429709)
\curveto(531.23884296,78.85429709)(535.94026296,101.47290709)(551.88861296,101.18646709)
\curveto(567.83695296,100.90002709)(595.36274296,96.11047709)(590.78607296,114.60055709)
\curveto(586.20939296,133.09062709)(595.39935296,132.31197709)(599.26703296,142.67363709)
\curveto(603.13470296,153.03529709)(602.53232296,153.65971709)(602.53232296,153.65971709)
}
}
{
\newrgbcolor{curcolor}{0 0 0}
\pscustom[linewidth=1.875,linecolor=curcolor,linestyle=dashed,dash=3.75 3.75]
{
\newpath
\moveto(531.35572296,77.59711709)
\lineto(474.36984296,68.77054709)
\lineto(474.36984296,68.77054709)
\lineto(474.36984296,68.77054709)
\lineto(474.36984296,68.77054709)
}
}
{
\newrgbcolor{curcolor}{1 1 1}
\pscustom[linestyle=none,fillstyle=solid,fillcolor=curcolor]
{
\newpath
\moveto(114.15552783,223.79080891)
\curveto(114.11963164,221.79219606)(112.43580786,220.19985262)(110.3969926,220.23647086)
\curveto(108.35817735,220.27308911)(106.73261984,221.92487087)(106.76851603,223.92348372)
\curveto(106.80441222,225.92209657)(108.488236,227.51444001)(110.52705126,227.47782176)
\curveto(112.56586652,227.44120352)(114.19142402,225.78942176)(114.15552783,223.79080891)
\closepath
}
}
{
\newrgbcolor{curcolor}{0 0 0}
\pscustom[linewidth=1.87500004,linecolor=curcolor]
{
\newpath
\moveto(114.15552783,223.79080891)
\curveto(114.11963164,221.79219606)(112.43580786,220.19985262)(110.3969926,220.23647086)
\curveto(108.35817735,220.27308911)(106.73261984,221.92487087)(106.76851603,223.92348372)
\curveto(106.80441222,225.92209657)(108.488236,227.51444001)(110.52705126,227.47782176)
\curveto(112.56586652,227.44120352)(114.19142402,225.78942176)(114.15552783,223.79080891)
\closepath
}
}
{
\newrgbcolor{curcolor}{0 0 0}
\pscustom[linewidth=4.375,linecolor=curcolor]
{
\newpath
\moveto(161.52829596,231.53031709)
\curveto(161.52829596,231.53031709)(178.21122996,237.98024709)(180.98076096,221.36340709)
\curveto(183.75029296,204.74657709)(184.77878396,193.68338709)(196.73080496,210.03580709)
\curveto(208.68282296,226.38822709)(180.43701696,259.41611709)(198.23657596,259.71002709)
\curveto(216.03613296,260.00392709)(227.29409296,237.71227709)(227.29409296,237.71227709)
}
}
{
\newrgbcolor{curcolor}{1 1 1}
\pscustom[linestyle=none,fillstyle=solid,fillcolor=curcolor]
{
\newpath
\moveto(163.99500583,231.48600591)
\curveto(163.95910964,229.48739306)(162.27528586,227.89504962)(160.2364706,227.93166786)
\curveto(158.19765535,227.96828611)(156.57209784,229.62006787)(156.60799403,231.61868072)
\curveto(156.64389022,233.61729357)(158.327714,235.20963701)(160.36652926,235.17301876)
\curveto(162.40534452,235.13640052)(164.03090202,233.48461876)(163.99500583,231.48600591)
\closepath
}
}
{
\newrgbcolor{curcolor}{0 0 0}
\pscustom[linewidth=1.87500004,linecolor=curcolor]
{
\newpath
\moveto(163.99500583,231.48600591)
\curveto(163.95910964,229.48739306)(162.27528586,227.89504962)(160.2364706,227.93166786)
\curveto(158.19765535,227.96828611)(156.57209784,229.62006787)(156.60799403,231.61868072)
\curveto(156.64389022,233.61729357)(158.327714,235.20963701)(160.36652926,235.17301876)
\curveto(162.40534452,235.13640052)(164.03090202,233.48461876)(163.99500583,231.48600591)
\closepath
}
}
{
\newrgbcolor{curcolor}{1 1 1}
\pscustom[linestyle=none,fillstyle=solid,fillcolor=curcolor]
{
\newpath
\moveto(230.98759183,237.64593291)
\curveto(230.95169564,235.64732006)(229.26787186,234.05497662)(227.2290566,234.09159486)
\curveto(225.19024135,234.12821311)(223.56468384,235.77999487)(223.60058003,237.77860772)
\curveto(223.63647622,239.77722057)(225.3203,241.36956401)(227.35911526,241.33294576)
\curveto(229.39793052,241.29632752)(231.02348802,239.64454576)(230.98759183,237.64593291)
\closepath
}
}
{
\newrgbcolor{curcolor}{0 0 0}
\pscustom[linewidth=1.87500004,linecolor=curcolor]
{
\newpath
\moveto(230.98759183,237.64593291)
\curveto(230.95169564,235.64732006)(229.26787186,234.05497662)(227.2290566,234.09159486)
\curveto(225.19024135,234.12821311)(223.56468384,235.77999487)(223.60058003,237.77860772)
\curveto(223.63647622,239.77722057)(225.3203,241.36956401)(227.35911526,241.33294576)
\curveto(229.39793052,241.29632752)(231.02348802,239.64454576)(230.98759183,237.64593291)
\closepath
}
}
{
\newrgbcolor{curcolor}{1 1 1}
\pscustom[linestyle=none,fillstyle=solid,fillcolor=curcolor]
{
\newpath
\moveto(266.70789683,244.98112261)
\curveto(266.67200064,242.98250976)(264.98817686,241.39016632)(262.9493616,241.42678456)
\curveto(260.91054635,241.46340281)(259.28498884,243.11518457)(259.32088503,245.11379742)
\curveto(259.35678122,247.11241027)(261.040605,248.70475371)(263.07942026,248.66813546)
\curveto(265.11823552,248.63151722)(266.74379302,246.97973546)(266.70789683,244.98112261)
\closepath
}
}
{
\newrgbcolor{curcolor}{0 0 0}
\pscustom[linewidth=1.87500004,linecolor=curcolor]
{
\newpath
\moveto(266.70789683,244.98112261)
\curveto(266.67200064,242.98250976)(264.98817686,241.39016632)(262.9493616,241.42678456)
\curveto(260.91054635,241.46340281)(259.28498884,243.11518457)(259.32088503,245.11379742)
\curveto(259.35678122,247.11241027)(261.040605,248.70475371)(263.07942026,248.66813546)
\curveto(265.11823552,248.63151722)(266.74379302,246.97973546)(266.70789683,244.98112261)
\closepath
}
}
{
\newrgbcolor{curcolor}{1 1 1}
\pscustom[linestyle=none,fillstyle=solid,fillcolor=curcolor]
{
\newpath
\moveto(412.30427783,49.69702691)
\curveto(412.26838164,47.69841406)(410.58455786,46.10607062)(408.5457426,46.14268886)
\curveto(406.50692735,46.17930711)(404.88136984,47.83108887)(404.91726603,49.82970172)
\curveto(404.95316222,51.82831457)(406.636986,53.42065801)(408.67580126,53.38403976)
\curveto(410.71461652,53.34742152)(412.34017402,51.69563976)(412.30427783,49.69702691)
\closepath
}
}
{
\newrgbcolor{curcolor}{0 0 0}
\pscustom[linewidth=1.87500004,linecolor=curcolor]
{
\newpath
\moveto(412.30427783,49.69702691)
\curveto(412.26838164,47.69841406)(410.58455786,46.10607062)(408.5457426,46.14268886)
\curveto(406.50692735,46.17930711)(404.88136984,47.83108887)(404.91726603,49.82970172)
\curveto(404.95316222,51.82831457)(406.636986,53.42065801)(408.67580126,53.38403976)
\curveto(410.71461652,53.34742152)(412.34017402,51.69563976)(412.30427783,49.69702691)
\closepath
}
}
{
\newrgbcolor{curcolor}{1 1 1}
\pscustom[linestyle=none,fillstyle=solid,fillcolor=curcolor]
{
\newpath
\moveto(478.92584783,69.36272691)
\curveto(478.88995164,67.36411406)(477.20612786,65.77177062)(475.1673126,65.80838886)
\curveto(473.12849735,65.84500711)(471.50293984,67.49678887)(471.53883603,69.49540172)
\curveto(471.57473222,71.49401457)(473.258556,73.08635801)(475.29737126,73.04973976)
\curveto(477.33618652,73.01312152)(478.96174402,71.36133976)(478.92584783,69.36272691)
\closepath
}
}
{
\newrgbcolor{curcolor}{0 0 0}
\pscustom[linewidth=1.87500004,linecolor=curcolor]
{
\newpath
\moveto(478.92584783,69.36272691)
\curveto(478.88995164,67.36411406)(477.20612786,65.77177062)(475.1673126,65.80838886)
\curveto(473.12849735,65.84500711)(471.50293984,67.49678887)(471.53883603,69.49540172)
\curveto(471.57473222,71.49401457)(473.258556,73.08635801)(475.29737126,73.04973976)
\curveto(477.33618652,73.01312152)(478.96174402,71.36133976)(478.92584783,69.36272691)
\closepath
}
}
{
\newrgbcolor{curcolor}{1 1 1}
\pscustom[linestyle=none,fillstyle=solid,fillcolor=curcolor]
{
\newpath
\moveto(534.29692783,77.57216691)
\curveto(534.26103164,75.57355406)(532.57720786,73.98121062)(530.5383926,74.01782886)
\curveto(528.49957735,74.05444711)(526.87401984,75.70622887)(526.90991603,77.70484172)
\curveto(526.94581222,79.70345457)(528.629636,81.29579801)(530.66845126,81.25917976)
\curveto(532.70726652,81.22256152)(534.33282402,79.57077976)(534.29692783,77.57216691)
\closepath
}
}
{
\newrgbcolor{curcolor}{0 0 0}
\pscustom[linewidth=1.87500004,linecolor=curcolor]
{
\newpath
\moveto(534.29692783,77.57216691)
\curveto(534.26103164,75.57355406)(532.57720786,73.98121062)(530.5383926,74.01782886)
\curveto(528.49957735,74.05444711)(526.87401984,75.70622887)(526.90991603,77.70484172)
\curveto(526.94581222,79.70345457)(528.629636,81.29579801)(530.66845126,81.25917976)
\curveto(532.70726652,81.22256152)(534.33282402,79.57077976)(534.29692783,77.57216691)
\closepath
}
}
{
\newrgbcolor{curcolor}{0 0 0}
\pscustom[linewidth=4.375,linecolor=curcolor]
{
\newpath
\moveto(460.25312296,226.16505709)
\curveto(487.25366296,226.29371709)(509.38005296,228.35069709)(507.84478296,211.19757709)
\curveto(506.30951296,194.04447709)(515.96964296,185.28062709)(523.29736296,183.30822709)
\curveto(530.62508296,181.33583709)(532.63800296,156.75584709)(532.63800296,156.75584709)
}
}
{
\newrgbcolor{curcolor}{1 1 1}
\pscustom[linestyle=none,fillstyle=solid,fillcolor=curcolor]
{
\newpath
\moveto(536.93388783,156.06507991)
\curveto(536.89799164,154.06646706)(535.21416786,152.47412362)(533.1753526,152.51074186)
\curveto(531.13653735,152.54736011)(529.51097984,154.19914187)(529.54687603,156.19775472)
\curveto(529.58277222,158.19636757)(531.266596,159.78871101)(533.30541126,159.75209276)
\curveto(535.34422652,159.71547452)(536.96978402,158.06369276)(536.93388783,156.06507991)
\closepath
}
}
{
\newrgbcolor{curcolor}{0 0 0}
\pscustom[linewidth=1.87500004,linecolor=curcolor]
{
\newpath
\moveto(536.93388783,156.06507991)
\curveto(536.89799164,154.06646706)(535.21416786,152.47412362)(533.1753526,152.51074186)
\curveto(531.13653735,152.54736011)(529.51097984,154.19914187)(529.54687603,156.19775472)
\curveto(529.58277222,158.19636757)(531.266596,159.78871101)(533.30541126,159.75209276)
\curveto(535.34422652,159.71547452)(536.96978402,158.06369276)(536.93388783,156.06507991)
\closepath
}
}
{
\newrgbcolor{curcolor}{1 1 1}
\pscustom[linestyle=none,fillstyle=solid,fillcolor=curcolor]
{
\newpath
\moveto(606.82820783,152.96894691)
\curveto(606.79231164,150.97033406)(605.10848786,149.37799062)(603.0696726,149.41460886)
\curveto(601.03085735,149.45122711)(599.40529984,151.10300887)(599.44119603,153.10162172)
\curveto(599.47709222,155.10023457)(601.160916,156.69257801)(603.19973126,156.65595976)
\curveto(605.23854652,156.61934152)(606.86410402,154.96755976)(606.82820783,152.96894691)
\closepath
}
}
{
\newrgbcolor{curcolor}{0 0 0}
\pscustom[linewidth=1.87500004,linecolor=curcolor]
{
\newpath
\moveto(606.82820783,152.96894691)
\curveto(606.79231164,150.97033406)(605.10848786,149.37799062)(603.0696726,149.41460886)
\curveto(601.03085735,149.45122711)(599.40529984,151.10300887)(599.44119603,153.10162172)
\curveto(599.47709222,155.10023457)(601.160916,156.69257801)(603.19973126,156.65595976)
\curveto(605.23854652,156.61934152)(606.86410402,154.96755976)(606.82820783,152.96894691)
\closepath
}
}
{
\newrgbcolor{curcolor}{1 1 1}
\pscustom[linestyle=none,fillstyle=solid,fillcolor=curcolor]
{
\newpath
\moveto(385.21489783,249.60222651)
\curveto(385.17900164,247.60361366)(383.49517786,246.01127022)(381.4563626,246.04788846)
\curveto(379.41754735,246.08450671)(377.79198984,247.73628847)(377.82788603,249.73490132)
\curveto(377.86378222,251.73351417)(379.547606,253.32585761)(381.58642126,253.28923936)
\curveto(383.62523652,253.25262112)(385.25079402,251.60083936)(385.21489783,249.60222651)
\closepath
}
}
{
\newrgbcolor{curcolor}{0 0 0}
\pscustom[linewidth=1.87500004,linecolor=curcolor]
{
\newpath
\moveto(385.21489783,249.60222651)
\curveto(385.17900164,247.60361366)(383.49517786,246.01127022)(381.4563626,246.04788846)
\curveto(379.41754735,246.08450671)(377.79198984,247.73628847)(377.82788603,249.73490132)
\curveto(377.86378222,251.73351417)(379.547606,253.32585761)(381.58642126,253.28923936)
\curveto(383.62523652,253.25262112)(385.25079402,251.60083936)(385.21489783,249.60222651)
\closepath
}
}
{
\newrgbcolor{curcolor}{1 1 1}
\pscustom[linestyle=none,fillstyle=solid,fillcolor=curcolor]
{
\newpath
\moveto(463.31119783,224.88292991)
\curveto(463.27530164,222.88431706)(461.59147786,221.29197362)(459.5526626,221.32859186)
\curveto(457.51384735,221.36521011)(455.88828984,223.01699187)(455.92418603,225.01560472)
\curveto(455.96008222,227.01421757)(457.643906,228.60656101)(459.68272126,228.56994276)
\curveto(461.72153652,228.53332452)(463.34709402,226.88154276)(463.31119783,224.88292991)
\closepath
}
}
{
\newrgbcolor{curcolor}{0 0 0}
\pscustom[linewidth=1.87500004,linecolor=curcolor]
{
\newpath
\moveto(463.31119783,224.88292991)
\curveto(463.27530164,222.88431706)(461.59147786,221.29197362)(459.5526626,221.32859186)
\curveto(457.51384735,221.36521011)(455.88828984,223.01699187)(455.92418603,225.01560472)
\curveto(455.96008222,227.01421757)(457.643906,228.60656101)(459.68272126,228.56994276)
\curveto(461.72153652,228.53332452)(463.34709402,226.88154276)(463.31119783,224.88292991)
\closepath
}
}
{
\newrgbcolor{curcolor}{1 1 1}
\pscustom[linestyle=none,fillstyle=solid,fillcolor=curcolor]
{
\newpath
\moveto(422.56607783,245.24980351)
\curveto(422.53018164,243.25119066)(420.84635786,241.65884722)(418.8075426,241.69546546)
\curveto(416.76872735,241.73208371)(415.14316984,243.38386547)(415.17906603,245.38247832)
\curveto(415.21496222,247.38109117)(416.898786,248.97343461)(418.93760126,248.93681636)
\curveto(420.97641652,248.90019812)(422.60197402,247.24841636)(422.56607783,245.24980351)
\closepath
}
}
{
\newrgbcolor{curcolor}{0 0 0}
\pscustom[linewidth=1.87500004,linecolor=curcolor]
{
\newpath
\moveto(422.56607783,245.24980351)
\curveto(422.53018164,243.25119066)(420.84635786,241.65884722)(418.8075426,241.69546546)
\curveto(416.76872735,241.73208371)(415.14316984,243.38386547)(415.17906603,245.38247832)
\curveto(415.21496222,247.38109117)(416.898786,248.97343461)(418.93760126,248.93681636)
\curveto(420.97641652,248.90019812)(422.60197402,247.24841636)(422.56607783,245.24980351)
\closepath
}
}
{
\newrgbcolor{curcolor}{1 1 1}
\pscustom[linestyle=none,fillstyle=solid,fillcolor=curcolor]
{
\newpath
\moveto(100.47028783,76.77358691)
\curveto(100.43439164,74.77497406)(98.75056786,73.18263062)(96.7117526,73.21924886)
\curveto(94.67293735,73.25586711)(93.04737984,74.90764887)(93.08327603,76.90626172)
\curveto(93.11917222,78.90487457)(94.802996,80.49721801)(96.84181126,80.46059976)
\curveto(98.88062652,80.42398152)(100.50618402,78.77219976)(100.47028783,76.77358691)
\closepath
}
}
{
\newrgbcolor{curcolor}{0 0 0}
\pscustom[linewidth=1.87500004,linecolor=curcolor]
{
\newpath
\moveto(100.47028783,76.77358691)
\curveto(100.43439164,74.77497406)(98.75056786,73.18263062)(96.7117526,73.21924886)
\curveto(94.67293735,73.25586711)(93.04737984,74.90764887)(93.08327603,76.90626172)
\curveto(93.11917222,78.90487457)(94.802996,80.49721801)(96.84181126,80.46059976)
\curveto(98.88062652,80.42398152)(100.50618402,78.77219976)(100.47028783,76.77358691)
\closepath
}
}
{
\newrgbcolor{curcolor}{1 1 1}
\pscustom[linestyle=none,fillstyle=solid,fillcolor=curcolor]
{
\newpath
\moveto(125.32215783,59.76015691)
\curveto(125.28626164,57.76154406)(123.60243786,56.16920062)(121.5636226,56.20581886)
\curveto(119.52480735,56.24243711)(117.89924984,57.89421887)(117.93514603,59.89283172)
\curveto(117.97104222,61.89144457)(119.654866,63.48378801)(121.69368126,63.44716976)
\curveto(123.73249652,63.41055152)(125.35805402,61.75876976)(125.32215783,59.76015691)
\closepath
}
}
{
\newrgbcolor{curcolor}{0 0 0}
\pscustom[linewidth=1.87500004,linecolor=curcolor]
{
\newpath
\moveto(125.32215783,59.76015691)
\curveto(125.28626164,57.76154406)(123.60243786,56.16920062)(121.5636226,56.20581886)
\curveto(119.52480735,56.24243711)(117.89924984,57.89421887)(117.93514603,59.89283172)
\curveto(117.97104222,61.89144457)(119.654866,63.48378801)(121.69368126,63.44716976)
\curveto(123.73249652,63.41055152)(125.35805402,61.75876976)(125.32215783,59.76015691)
\closepath
}
}
{
\newrgbcolor{curcolor}{1 1 1}
\pscustom[linestyle=none,fillstyle=solid,fillcolor=curcolor]
{
\newpath
\moveto(163.20961183,51.10292691)
\curveto(163.17371564,49.10431406)(161.48989186,47.51197062)(159.4510766,47.54858886)
\curveto(157.41226135,47.58520711)(155.78670384,49.23698887)(155.82260003,51.23560172)
\curveto(155.85849622,53.23421457)(157.54232,54.82655801)(159.58113526,54.78993976)
\curveto(161.61995052,54.75332152)(163.24550802,53.10153976)(163.20961183,51.10292691)
\closepath
}
}
{
\newrgbcolor{curcolor}{0 0 0}
\pscustom[linewidth=1.87500004,linecolor=curcolor]
{
\newpath
\moveto(163.20961183,51.10292691)
\curveto(163.17371564,49.10431406)(161.48989186,47.51197062)(159.4510766,47.54858886)
\curveto(157.41226135,47.58520711)(155.78670384,49.23698887)(155.82260003,51.23560172)
\curveto(155.85849622,53.23421457)(157.54232,54.82655801)(159.58113526,54.78993976)
\curveto(161.61995052,54.75332152)(163.24550802,53.10153976)(163.20961183,51.10292691)
\closepath
}
}
{
\newrgbcolor{curcolor}{1 1 1}
\pscustom[linestyle=none,fillstyle=solid,fillcolor=curcolor]
{
\newpath
\moveto(204.34032683,52.20498691)
\curveto(204.30443064,50.20637406)(202.62060686,48.61403062)(200.5817916,48.65064886)
\curveto(198.54297635,48.68726711)(196.91741884,50.33904887)(196.95331503,52.33766172)
\curveto(196.98921122,54.33627457)(198.673035,55.92861801)(200.71185026,55.89199976)
\curveto(202.75066552,55.85538152)(204.37622302,54.20359976)(204.34032683,52.20498691)
\closepath
}
}
{
\newrgbcolor{curcolor}{0 0 0}
\pscustom[linewidth=1.87500004,linecolor=curcolor]
{
\newpath
\moveto(204.34032683,52.20498691)
\curveto(204.30443064,50.20637406)(202.62060686,48.61403062)(200.5817916,48.65064886)
\curveto(198.54297635,48.68726711)(196.91741884,50.33904887)(196.95331503,52.33766172)
\curveto(196.98921122,54.33627457)(198.673035,55.92861801)(200.71185026,55.89199976)
\curveto(202.75066552,55.85538152)(204.37622302,54.20359976)(204.34032683,52.20498691)
\closepath
}
}
{
\newrgbcolor{curcolor}{1 1 1}
\pscustom[linestyle=none,fillstyle=solid,fillcolor=curcolor]
{
\newpath
\moveto(587.60485783,178.64091291)
\curveto(587.56896164,176.64230006)(585.88513786,175.04995662)(583.8463226,175.08657486)
\curveto(581.80750735,175.12319311)(580.18194984,176.77497487)(580.21784603,178.77358772)
\curveto(580.25374222,180.77220057)(581.937566,182.36454401)(583.97638126,182.32792576)
\curveto(586.01519652,182.29130752)(587.64075402,180.63952576)(587.60485783,178.64091291)
\closepath
}
}
{
\newrgbcolor{curcolor}{0 0 0}
\pscustom[linewidth=1.87500004,linecolor=curcolor]
{
\newpath
\moveto(587.60485783,178.64091291)
\curveto(587.56896164,176.64230006)(585.88513786,175.04995662)(583.8463226,175.08657486)
\curveto(581.80750735,175.12319311)(580.18194984,176.77497487)(580.21784603,178.77358772)
\curveto(580.25374222,180.77220057)(581.937566,182.36454401)(583.97638126,182.32792576)
\curveto(586.01519652,182.29130752)(587.64075402,180.63952576)(587.60485783,178.64091291)
\closepath
}
}
{
\newrgbcolor{curcolor}{0 0 0}
\pscustom[linewidth=1.87499991,linecolor=curcolor]
{
\newpath
\moveto(237.36740227,185.49465846)
\curveto(236.91686436,185.02548343)(225.84806515,168.45679221)(221.60754923,162.58117477)
\curveto(220.27196213,160.73061047)(218.47937609,159.3235743)(216.9153046,157.69477879)
}
}
{
\newrgbcolor{curcolor}{0 0 0}
\pscustom[linewidth=1.87499991,linecolor=curcolor]
{
\newpath
\moveto(233.95990591,200.30521662)
\curveto(232.77314604,198.66087298)(242.45522979,174.64955473)(245.66935158,168.61389366)
\curveto(247.24362773,165.65762678)(246.39942212,165.18209008)(247.08469764,161.96090487)
}
}

 \uput[0](34,146){$e$}
\uput[0](284,140){$g$}
\uput[0](26,268){$\mc H$}
\uput[0](594,268){$\mc K$}

\end{pspicture}

%% file: sos.bbl
\begin{thebibliography}{10}

\bibitem{wagner-choe}
Y.~B. Choe and D.~G. Wagner.
\newblock Rayleigh matroids.
\newblock {\em Comb. Probab. Comput.}, 15(5):765--781, 2006.

\bibitem{cocks}
C.~C. Cocks.
\newblock Correlated matroids.
\newblock {\em Combin. Probab. Comput.}, 17(4):511--518, 2008.

\bibitem{Erickson-2008}
A.~Erickson.
\newblock Negative correlation properties for matroids.
\newblock Master's thesis, University of Waterloo, 2008,
  \url{http://hdl.handle.net/10012/4165}.

\bibitem{fkg}
C.~M. Fortuin, J.~Ginibre, and P.~W. Kasteleyn.
\newblock Correlation inequalities on some partially ordered sets.
\newblock {\em Comm. Math. Phys.}, 22(2):89--103, 1971.

\bibitem{grimmett-winkler}
G.~R. Grimmett and S.~N. Winkler.
\newblock Negative association in uniform forests and connected graphs.
\newblock {\em Random Structures Algorithms}, 24(4):444--460, 2004.

\bibitem{kahnforests}
J.~Kahn.
\newblock A normal law for matchings.
\newblock {\em Combinatorica}, 20(3):339--391, 2000.

\bibitem{kahn-neiman}
J.~Kahn and M.~Neiman.
\newblock {Negative correlation and log-concavity}.
\newblock {\em Random Structures and Algorithms (to appear)}, 2009.

\bibitem{kirchhoff}
G.~Kirchhoff.
\newblock Uber die aufl\"osung der gleichungen, auf welche man bei der
  untersuchungen der linearen vertheilung galvanischer str\"ome gef\"uhrt wird.
\newblock {\em Ann. Phys. Chem.}, 72:497--508, 1847.

\bibitem{merino}
C.~Merino.
\newblock {\em Matroids, the Tutte polynomial and the chip firing game}.
\newblock PhD thesis, Somerville College, University of Oxford, 1999.

\bibitem{pemantle}
R.~Pemantle.
\newblock Towards a theory of negative dependence.
\newblock {\em J. Math. Phys.}, 41:1371--1390, 2000.

\bibitem{semple-welsh}
C.~Semple and D.~Welsh.
\newblock Negative correlation in graphs and matroids.
\newblock {\em Comb. Probab. Comput.}, 17(3):423--435, 2008.

\bibitem{seymour-regular}
P.~D. Seymour.
\newblock Decomposition of regular matroids.
\newblock {\em Journal of Combinatorial Theory Series B}, 28:305--359, 1980.

\bibitem{sflows}
P.~D. Seymour.
\newblock Matroids and multicommodity flows.
\newblock {\em European J. Combin.}, 2(3):257--290, 1981.

\bibitem{ncrv}
D.~G. Wagner.
\newblock Negatively correlated random variables and {M}ason's conjecture for
  independent sets in matroids.
\newblock {\em Ann. Comb.}, 12(2):211--239, 2008.

\end{thebibliography}
